\documentclass{emsprocart}

\usepackage{amsmath}
\usepackage{amssymb}
\usepackage[all]{xy}

\contact[manin@mpim-bonn.mpg.de]{Yuri I.~Manin, Max--Planck--Institut f\"ur Mathematik, Vivatsgasse 7, Bonn, D-53111, Germany}
\contact[matilde@caltech.edu]{Matilde Marcolli, Department of Mathematics,
California Institute of Technology, 1200 E California Boulevard, Pasadena, CA 91125, USA}

\newtheorem{theorem}{Theorem}[section]

\newtheorem{lemma}[theorem]{Lemma}
\newtheorem{proposition}[theorem]{Proposition}

\theoremstyle{definition}
\newtheorem{definition}[theorem]{Definition}
\newtheorem{remark}[theorem]{Remark}

\newtheorem{example}[theorem]{Example}

\def\F{{\mathbb F}}
\def\Z{{\mathbb Z}}
\def\cV{{\mathcal V}}
\def\bT{{\mathbb T}}
\def\bG{{\mathbb G}}
\def\A{{\mathbb A}}
\def\bL{{\mathbb L}}
\def\roman{\rm}
\renewcommand{\P}{{\mathbb P}}
\def\Q{{\mathbb Q}}
\def\cT{{\mathcal T}}
\def\cC{{\mathcal C}}
\def\C{{\mathbb C}}
\def\cG{{\mathcal G}}

\title[Moduli Operad over $\F_1$]{Moduli Operad over $\F_1$}

\author[Yuri I.~Manin and Matilde Marcolli]{Yuri I.~Manin and Matilde Marcolli
\thanks{The second author acknowledges support and hospitality of the 
Max Planck Institute and the Mathematical Sciences Research Institute and
support from NSF grants DMS-0901221, DMS-1007207, DMS-1201512, 
PHY-1205440. We thank Paolo Aluffi, Tom Graber and Oliver Lorscheid for
constructive criticism, useful comments and discussions.}}

\begin{document}

\begin{abstract}
In this paper we answer a question raised in \cite{Man1}, 
Sec.~4, by showing that the genus zero moduli operad $\{\overline{M}_{0,n+1}\}$
can be endowed with natural descent data that allow it to be considered as  the lift to 
${\rm Spec}\,\Z$ of an operad over $\F_1$. The relevant descent data are based
on a notion of constructible sets and constructible functions over $\F_1$, which
describes suitable differences of torifications with a positivity condition on the
class in the Grothendieck ring. More generally, we do the same for the 
operads $\{T_{d,n+1}\}$ (whose components were) introduced in \cite{CGK}.
Finally, we describe a blueprint structure on $\{\overline{M}_{0,n}\}$ and we
discuss from this perspective the genus zero boundary modular operad 
$\{\overline{M}_{g,n+1}^0\}$.
\end{abstract}

\tableofcontents

\section{Introduction and summary}

Of many recently suggested definitions
of $\F_1$--geometry, we work  with the one developed in \cite{LL1} that seems to be the
minimal one. Namely, an $\F_1$--scheme is represented by its lift to ${\rm Spec}(\Z)$
and  the relevant descent data which are essentially a representation of the lifted scheme
as a disjoint union of locally closed tori. 

\smallskip

This notion of $\F_1$--geometry can be seen as the simplest geometrization of
the condition that the class in the Grothendieck ring of the variety decomposes
as a sum of classes of tori, with non-negative coefficients. This motivic condition
accounts for the expected behavior of points over $\F_1$ and over ``extensions"
$\F_{1^m}$ in relation to the counting of points over $\F_{q^m}$ and zeta functions.

\smallskip

In this setting, we show that, while the torification condition (possibly with additional
restrictions such as a compatibility with an affine covering) provides a viable notion of
``algebraic variety over $\F_1$", when one considers possible descent data to $\F_1$
for stable curves of genus zero with marked points one needs to consider also
objects that are analogs of ``constructible sets" over $\F_1$, which can be seen as 
formal differences of torifications. In general, the complement of an algebraic variety
in another need not be an algebraic variety, but it is a constructible set. Similarly,
not all points or subvarieties over $\F_1$ (in the sense of torifications as well as
in other forms of $\F_1$-geometry) are {\em complemented}. The complemented
case corresponds to those $\F_1$-subarieties whose complement also defines
an $\F_1$-variety, while in the non-complemented case one obtains an $\F_1$-constructible
set, according to a suitable notion of differences of torifications that we refer to
as ``constructible torifications". The moduli spaces $\overline{M}_{0,n}$ and
their generalizations $T_{d,n}$ constructed in \cite{CGK} also have a structure
of $\F_1$-constructible sets. The operad structure on these moduli spaces
is also compatible with the $\F_1$-structure and the operad morphisms give
rise to $\F_1$-constructible morphisms.

\smallskip

In Section \ref{torSec}, we recall the notion of torification from \cite{LL1} and
we discuss different equivalence relations that determine when two choices
of torification on the same variety over $\Z$ determine the same $\F_1$-structure.
This leads to three different notions of $\F_1$-morphisms, which we refer to
as strong, ordinary, and weak morphisms. 

\smallskip

In Section \ref{GrSec} we focus on the condition that the Grothendieck
class of a variety decomposes into a sum of tori with non-negative coefficients,
which is necessary for the existence of geometric torifications.
We show that it is satisfied for the moduli spaces $\overline{M}_{0,n}$ and
$T_{d,n}$. This follows the same argument used in \cite{Man2} and \cite{CGK},
respectively, for the computation of the Poincar\'e polynomials. We also show
that these same computations provide a generating series for the numbers 
of $F_{1^m}$--points of $\overline{M}_{0,n}$ and $T_{d,n}$.

\smallskip

In Sections \ref{ConstrSec} and \ref{ConstrSec2} we discuss the notion
of {\em complemented} points and {\em complemented subspaces}
in $\F_1$-geometry. We analyze the geometric torifications of stable curves of
genus zero and the role of the marked points as uncomplemented
points. We introduce a notion of {\em constructible sets} over $\F_1$ and
of {\em constructible torifications}, which are formal differences of
torifications preserving the positivity of the Grothendieck class. In Section
\ref{ctorSec} we show that the moduli spaces $\overline{M}_{0,n}$ and
$T_{d,n}$ are $\F_1$-constructible sets.

\smallskip

In Section \ref{opSec}, for each $d\ge 1$, we introduce the operads  
with components $\{T_{d,n+1}\}$ from \cite{CGK} and we show that
the operadic structure morphisms are compatible with the structure
of $\F_1$-constructible sets. The operad composition operations
and the morphisms that forget marked points determine strong
$\F_1$-constructible morphisms, while the action of $S_n$ that
permutes marked points acts through ordinary $\F_1$-constructible 
morphisms. In Section \ref{wondSec} we also show that, if one
uses the description of the moduli spaces $\overline{M}_{0,n}$ and
$T_{d,n}$ as iterated blowups, related to the Fulton--MacPherson
compactifications as in \cite{CGK}, then the projection maps of
the iterated blowups are only weak $\F_1$-morphisms.

\smallskip

In Section \ref{blueSec} we focus on the {\em blueprint} approach to
$\F_1$-geometry, developed in \cite{Lo}, \cite{Lo2}. We 
make explicit a blueprint structure of   $\overline{M}_{0,n}$ based upon
explicit equations for $\overline{M}_{0,n}$, as in \cite{GiM1}, \cite{KeTe}. 
We consider then 
the genus zero boundary modular operad $\{\overline{M}_{g,n+1}^0\}$ whose 
components are, by definition, unions of those boundary strata in 
$\{\overline{M}_{g,n+1}\}$ that parametrize curves whose normalized 
irreducible components are projective lines. 
This is an operad in the category of DM--stacks, so that for its complete treatment
within the setting of torifications 
it would be necessary to develop a formalism of stacky $\F_1$--geometry
compatible with torifications as descent data. 
We describe a blueprint structure on the genus zero boundary 
$\overline{M}_{g,n+1}^0$ of the higher genus moduli spaces, using
a crossed product construction. 

\smallskip

\bigskip

\section{Torifications}\label{torSec}

Of many recently suggested definitions
of $\F_1$--geometry, we work  with the one developed in \cite{LL1} that seems to be the
minimal one. Namely, an $\F_1$--scheme is represented by its lift to $\roman{Spec}\,\Z$
and  the relevant descent data which are essentially a representation of the lifted scheme
as a disjoint union of locally closed tori. 

\smallskip

The notion of torification introduced in  \cite{LL1} is the following condition, which
we refer to in this paper as {\em geometric torification}.

\begin{definition}\label{def1.1}
(\cite{LL1}).  
A torification of the scheme $X$  is a morphism of schemes $e_X : T \to X$ from a disjoint union of tori 
  $T=\coprod_{j\in I} T_i$, $T_j=\bG_m^{d_j}$,
such that the restriction of $e_X$  to each torus is an immersion (i.e.~isomorphism with a
locally closed subscheme), 
and $e_X$ induces bijections of $k$--points, $e_X(k): T(k) \to X(k)$, for every field $k$.
\end{definition}

Moreover, in  \cite{LL1} the authors also consider the stronger notion of
{\em affine torification}. 

\begin{definition}\label{def1.1aff}
(\cite{LL1}).  
The torification $e_X$ is called  affine if there exists an affine covering  $\{ U_\alpha \}$ of $X$ compatible with $e_X$ 
in the following sense: for each affine open set $U_\alpha$ in the covering, there is a 
subfamily of tori $\{T_j\,|\,j\in I_{\alpha}\}$ in the torification $e_X$  such that the restriction  of
$e_X$ to the disjoint union of tori from this subfamily  is a of torification of $U_\alpha$.
\end{definition}

\smallskip
\subsection{Levels of torified structures}\label{SecTor1}

We assume that $X$ is a variety over $\Z$.
There are three levels of increasingly restrictive conditions in this approach based on
defining $\F_1$-structures via torifications: the basic level is a decomposition of
the class in the Grothendieck ring, the second is a geometrization of this decomposition
at the level of the variety itself, and the third level includes more restrictive conditions,
such as affine and regular. 

\begin{enumerate}
\item {\em Torification of the Grothendieck class}: this is the weakest condition and it
simply consists of the requirement that the class $[X] \in K_0(\cV_\Z)$ in the Grothendieck
ring can be written as
\begin{equation}\label{K0torif}
[X] = \sum_k a_k \bT^k,
\end{equation}
where $\bT=[\bG_m]=\bL-1$, and with coefficients $a_k \geq 0$.
\item {\em Geometric torification}: this is the condition of Definition \ref{def1.1} above.
\item {\em Affine torification}, where the geometric torification is also affine in the
sense of Definition \ref{def1.1aff}.
\item {\em Regular torification}: this is a geometric torification where one also 
requires that the closure of each torus in the torification is itself a 
union of tori of the torification.
\end{enumerate}

\smallskip

Roughly, one can understand these different levels as describing stronger forms of
$\F_1$-structures based on torification. The decomposition of the class in the Grothendieck
ring reflects how one expects that $\Z$-varieties that descend to $\F_1$ should behave with
respect to motivic properties such as the zeta function and counting of points. The notion
of geometric torification introduced in  \cite{LL1} can be seen as a minimal way of
making this motivic behavior ``geometric". The further level, given by the affine condition,
was introduced in  \cite{LL1}, motivated by the comparison between this approach to
$\F_1$-geometry and the approaches developed by Soul\'e \cite{Soule} and Connes--Consani
\cite{CC1}. However, varieties like Grassmannians, which in many respects it would be
natural to expect should descend to $\F_1$, have natural torifications coming from their
cell decompositions that are not affine. This concern justifies retaining the intermediate
level of $\F_1$-structure given by geometric torifications {\em without} the affine condition.
As we shall argue later, this level already provides a very rich and interesting structure.
The regularity condition, which is independent of the affine requirement, but is usually
considered for affine torifications, was introduced in  \cite{LL1} as a possible way to ``rigidify" the
choice of torification. We follow here a different approach based on considering different
levels of equivalence relations among torifications, hence we will not consider the
regularity condition.

\medskip

\subsection{Equivalent torifications and morphisms}\label{eqtorSec}

When we consider {\em geometric torifications} as data defining $\F_1$-structures
on $\Z$-varieties, one would like to have a natural equivalence relation describing
when two different choices of torification on the same varieties should be
regarded as defining the {\em same} $\F_1$-structure.

\smallskip

We first recall the notion of {\em torified morphism} introduced in \cite{LL1}.

\begin{definition}\label{def1.1mor} (\cite{LL1}). 
A morphism of torified varieties (torified morphism) $\Phi:(X,e_X:\, T_X\to X)\to (Y,e_Y:\,T_Y\to Y)$ is a
triple $\Phi=(\phi,\psi,\{ \phi_i \})$ where $\phi: X \to Y$ is a morphism
of $\Z$-varieties, $\psi: I_X\to I_Y$ is a map of the indexing sets
of the two torifications, and $\phi_j: T_{X,j} \to T_{Y,\psi(j)}$ is a morphism of
algebraic groups, such that $\phi\circ e_X|_{T_{X,j}} = e_Y|_{T_{Y,\psi(j)}} \circ \phi_j$.
\end{definition}

In \cite{LL1}, a notion of {\em affinely torified morphism} was also introduced: these
are torified morphisms in the sense recalled above, between affinely torified varieties,
such that, if $\{ U_j \}$ is an affine open covering of $X$ compatible with the
torification, then for every $j$ the image of $U_j$ under $\Phi$ is 
an affine subscheme of $Y$. The following lemma, communicated to us by
Lorscheid, shows that it is not necessary to assume this as an additional 
condition for torified morphisms between affine affinely torified varieties.

\begin{lemma}\label{affmorphs}
Let $\Phi:(X,e_X:\, T_X\to X)\to (Y,e_Y:\,T_Y\to Y)$ be a torified morphism
between affinely torified varieties, with $\{ U_i \}$ and $\{ V_j \}$ respective 
affine torified coverings. Then $\Phi$ is an affinely torified morphism.
\end{lemma}

\proof
Let $W_j=\Phi^{-1}(V_j)$ and $U_{ij}=U_i\cap W_j$. Then $U_{ij}$ is a torified and
quasi-affine subscheme of $X$ that maps to $V_j$. The collection of all
$U_{ij}$ covers $X$. Consider then $Z_{ij}={\rm Spec}\, O_X(U_{ij})$. Since
$U_{ij}$ is quasi-affine, the natural map $U_{ij} \to Z_{ij}$ is an embedding 
of $U_{ij}$ into the affine subscheme $Z_{ij}$ of $U_i$. Moreover, the 
morphism $U_{ij} \to V_j$ extends naturally to a morphism $Z_{ij} \to V_j$ 
(since $V_j$ is affine), which means that $Z_{ij}$ is contained in $W_j$. 
Therefore, $Z_{ij}$ is contained in $U_i \cap W_j = U_{ij}$, hence 
we have that $U_{ij}=Z_{ij}$ is affine and therefore $\Phi$ is an affinely
torified morphism.
\endproof

\medskip

We consider the following notions of equivalence of torifications on a given $\Z$-variety $X$.

\begin{enumerate}
\item {\em Strong equivalence}: the identity morphism is torified.
\item {\em Ordinary Equivalence}: there exists an isomorphism of $X$ that is torified.
\item {\em Weak equivalence}: one identifies as the same $\F_1$-structure two
torifications on a variety $X$ such that $X$ has a decompositions into a disjoint union
of subvarieties  $X=\cup_j X_j$ and $X=\cup_j X'_j$, respectively compatible with the 
torifications, and such that there exist isomorphisms $\phi_i: X_i \to X'_i$ that are torified.
One considers the equivalence relation generated by these identifications.
\end{enumerate}

In the case of a weak equivalence the isomorphisms on the pieces of the decomposition
do not necessarily extend to isomorphisms of the whole variety.
Typical examples of this third condition are obtained by considering cell decompositions
compatible with the torifications. For example, one can consider $\P^1 \times \P^1$ with
the cell decomposition $\P^1 =\A^0 \cup \A^1$ on each factor. One can then consider
the standard torification of  $\P^1 \times \P^1$ compatible with the cell decomposition
and the torification obtained by taking a torification of the diagonal and of its complement 
in the $\A^2$ cell, and the torification of the other cells as before. These two torifications
are related by a weak equivalence but not by an ordinary one.

\medskip

The choice of the equivalence relation above determines what morphisms of
$\Z$-varieties can be regarded as descending to $\F_1$.

\begin{enumerate}
\item {\em Strong $\F_1$-morphisms} (or strongly torified morphisms): 
when geometric torifications are assumed to
define the same $\F_1$-structure iff they are strongly equivalent, morphisms of 
$\Z$-varieties that define $\F_1$-morphisms are torified morphisms in the sense of Definition
\ref{def1.1mor}.
\item {\em Ordinary $\F_1$-morphisms} (or ordinarily torified morphisms): under ordinary equivalence, then $\F_1$-morphisms
are all morphisms of $\Z$-varieties that become torified after composing with 
isomorphisms.
\item {\em Weak $\F_1$-morphisms} (or weakly torified morphisms): 
under weak equivalence, $\F_1$-morphisms
are morphisms of $\Z$-varieties that become torified after composition with
weak equivalences. 
\end{enumerate}

\medskip

In the following, we refer to the different cases above as a strong, ordinary, or weak
$\F_1$-structure or as geometric torifications in the strong, ordinary, or weak sense.

\medskip

 \begin{example}
 Any toric variety has a natural torification by torus orbits.  In \cite{LL1}, explicit
 affine torifications are constructed, and it is checked that  toric morphisms are compatible with them.
 This shows that the Losev--Manin operad $\{\overline{L}_{0,n}\}$ in \cite{LoMa}, \cite{LoMa2}, 
 \cite{Man1}   has natural descent data to $\F_1$, in the strong sense of the notion of 
 torifications and morphisms described above.
 \end{example}

\medskip

\begin{remark} 
Considering torifications and $\F_1$-morphisms in the weak sense is 
very close to imposing only the condition of torification of  
Grothendieck classes, though it appears to be stronger, as our
discussion of constructible torifications in \S \ref{ConstrSec} will illustrate.
\end{remark}

\medskip

\begin{remark}\label{bluermk}
Among the other existing approaches to $\F_1$-structures, the one based on
the notion of {\em blueprint}, developed in \cite{Lo}, \cite{Lo2}, \cite{LL3},
is not based on decompositions into tori, and it is a less restrictive form
of $\F_1$-structure in the sense that every scheme of finite type admits a
``blue model" of finite type over $\F_1$. 
\end{remark}

\medskip

\subsection{Categories of geometric torifications}

The different notions of morphisms of torified varieties considered
above lead to the following categorical formulation.

\begin{proposition}\label{catsTor}
There are categories $\cG\cT^s \subset \cG\cT^o \subset \cG\cT^w$ where
the objects, ${\rm Obj}(\cG\cT^s)={\rm Obj}(\cG\cT^o)={\rm Obj}(\cG\cT^w)$,
are pairs $(X_\Z, \cT)$, with $X_\Z$ a variety over $\Z$ and $\cT=\{ T_i \}$
a geometric torification of $X_\Z$.
Morphisms in $\cG\cT^s$ are strong morphisms of
geometrically torified spaces; morphisms in $\cG\cT^o$ are ordinary
morphisms of geometrically torified spaces; 
morphisms in $\cG\cT^w$ are weak 
morphisms of geometrically torified spaces.
\end{proposition}

\proof According to our previous discussion, strong morphisms
of geometrically torified spaces are the ``torified morphisms" of
Definition \ref{def1.1mor}, hence the category $\cG\cT^s$ is the
category of torified varieties, as considered in \cite{LL1}.
Morphisms in $\cG\cT^o$ are arbitrary compositions of  
torified morphisms and ordinary equivalences, which means
that they can be written as arbitrary compositions of 
torified morphisms and isomorphisms of $\Z$-varieties. 
Since composition of two such morphisms will still be of the
same kind, composition of morphisms is well defined in 
$\cG\cT^o$. Morphisms in $\cG\cT^s$ are also morphisms
in $\cG\cT^o$, but not the other way around.
Similarly, morphisms in $\cG\cT^w$ are arbitrary compositions
of torified morphisms and weak equivalences, that is, arbitrary compositions
of torified morphisms and local isomorphisms of the type described in
\S \ref{eqtorSec} above. Again, composition is well defined.
Morphisms in $\cG\cT^s$ and morphisms in $\cG\cT^o$ are also morphisms
in $\cG\cT^w$, but not conversely.
\endproof

\medskip

\section{Grothendieck classes and torifications}\label{GrSec}

In this section we consider the moduli spaces $\overline{M}_{0,n}$, as well
as their generalizations $T_{d,n}$ considered in \cite{CGK}, from the point
of view of classes in the Grothendieck ring. The existence of a decomposition
of the form \eqref{K0torif} into tori, with non-negative coefficients, follows from
the fact that these spaces can be realized as a sequence of iterated
blowups starting from a variety that clearly admits a torification and blowing
up loci that, in turn admit torifications. The explicit form of the decomposition
\eqref{K0torif} mirrors the known formulae for the Poincar\'e polynomial and
the Euler characteristic of \cite{Man2} and \cite{CGK} and can be obtained
by a similar argument. The generating functions of  \cite{Man2} and \cite{CGK}
computing the Poincar\'e polynomials are also related to counting points over
the extensions $\F_{1^m}$.

\medskip

\subsection{The class of $M_{0,n}$}

A first simple observation, which will be useful in the following, is that
the open stratum $M_{0,n}$ by itself cannot be torifed, since it fails the 
necessary condition that the class $[M_{0,n}]$ is torified by a 
decomposition \eqref{K0torif} with non-negative coefficients.

\begin{lemma}\label{notoriM0n}
The class $[ M_{0,n} ]$ has a decomposition into tori of the form
\begin{equation}\label{K0M0nTi}
 [ M_{0,n} ] = \sum_{k=0}^{n-2} s(n-2,k)\,  \sum_{j=0}^k \binom{k}{j} \bT^j,
\end{equation}
where $s(m,k)$ is the Stirling number of the first kind. In particular,
the open stratum $M_{0,n}$ does not admit a geometric torification. 
\end{lemma}

\proof We can view $M_{0,n}$ as the complement of the diagonals in a
product of $n-3$ copies of $\P^1\smallsetminus \{ 0,1,\infty \}$, hence
the class in the Grothendieck ring is given by
\begin{equation}\label{K0M0n}
[ M_{0,n} ] = (\bT-1) (\bT-2) \cdots (\bT-n+2)= \binom{\bT -1}{n-3} \, (n-3)! =(-1)^n (1-\bT)_{n-2}  ,
\end{equation}
where $(x)_m=\Gamma(x+m)/\Gamma(x)$ is the Pochhammer symbol, satisfying
$$(x)_m= \sum_{k=0}^m (-1)^{m-k} s(m,k)\, x^k,$$ with coefficients $s(m,k)$
the Stirling numbers of the first kind, so that $(-1)^{m-k} s(m,k)$ is the number
of permutations in $S_m$ consisting of $k$ cycles. 
Thus, we obtain 
$$ [ M_{0,n} ] = (-1)^n \sum_{k=0}^{n-2} (-1)^{n-k} s(n-2,k)\, 
(-1)^k (\bT-1)^k $$
which gives \eqref{K0M0nTi}, where some of the coefficients are clearly negative.
\endproof

\medskip

\subsection{The class of $\overline{M}_{0,n}$ and $\F_{1^m}$-points}

By Lemma \ref{notoriM0n}, the open stratum $M_{0,n}$ by itself 
cannot be torified. However, when one considers the compactification 
$\overline{M}_{0,n}$, one finds that the torification of the Grothendieck 
class is satisfied. 

\medskip

\begin{proposition}\label{K0barM0n}
The classes $[\overline{M}_{0,n}]\in K_0(\cV_\Z)$ fit into a generating series
\begin{equation}\label{K0M0nseries}
\varphi(t)= t + \sum_{n=2}^\infty [\overline{M}_{0,n}]\, \frac{t^n}{n!},
\end{equation}
in $K_0(\cV_\Z)_\Q[[t]]$, with $K_0(\cV_\Z)_\Q=K_0(\cV_\Z)\otimes_\Z \Q$,
where $\varphi(t)$ is the unique solution in
$t+t^2 K_0(\cV_\Z)_\Q[[t]]$ of the differential equation
\begin{equation}\label{K0eqdiff}
(1+\bL\, t - \bL\, \varphi(t)) \, \varphi^\prime(t) = 1 + \varphi(t).
\end{equation}
In particular, the classes $[\overline{M}_{0,n}]$ satisfy the recursive
relation
\begin{equation}\label{recK0M0n}
[\overline{M}_{0,n+2}]  = [\overline{M}_{0,n+1}] + \bL \, \sum_{i+j=n+1, i\geq 2} \binom{n}{i}
[\overline{M}_{0,i+1}]\, [\overline{M}_{0,j+1}] ,
\end{equation}
and therefore have a decomposition \eqref{K0torif} with non-negative coefficients.
\end{proposition}

\proof The argument is analogous to the proof of Theorem 0.3.1 of \cite{Man2} computing
the Poincar\'e polynomials of $\overline{M}_{0,n}$. In fact, the same argument used in
 \cite{Man2} to determine the Poincar\'e polynomials applies to the computation of the
 Grothendieck classes, using the classes of all the $M_{0,k}$ given in \eqref{K0M0n}, which
we rewrite as  
$$ [M_{0,k}] = \binom{\bL -2}{k-3} \, (k-3)! $$
which is the direct analog of equation (1.2) of \cite{Man2} for the Poincar\'e polynomials.
The existence of a decomposition \eqref{K0torif} with non-negative coefficients then 
follows inductively from the fact that the classes satisfy the recursive relation \eqref{recK0M0n}, 
which follows from \eqref{K0eqdiff} as in Corollary 0.3.2 of \cite{Man2}, and that the first 
terms of the recursion can be seen explicitly to have a non-negative coefficients.
\endproof

It would be interesting to know if the Chern class of $\overline{M}_{0,n}$
also satisfies a similar recursive formula and positivity property.

\begin{remark}\label{K0Poincare}
The Poincar\'e polynomial for $\overline{M}_{0,n}$ can be recovered from the Grothendieck
class by formally replacing $\bL$ with $q^2$ in the resulting expression. This fact holds
more generally for smooth projective varieties whose class in the Grothendieck ring is a
polynomial $[X]=\sum_k b_k \bL^k$ in the class $\bL$ of the Lefschetz motive. 
In fact, in this case the Hodge--Deligne polynomial $h_X(u,v)=\sum_{p,q} (-1)^{p+q} 
h^{p,q}(X_\C) u^p v^q$ is given by $h_X(u,v)=\sum_k b_k (uv)^k$, which implies that
$X_\C$ is Hodge-Tate, namely $h^{p,q}(X_\C)=0$ for $p\neq q$. This in turn implies
that the Poincar\'e polynomial is given by $P_X(q)=\sum_k b_k q^{2k}$, hence it is
obtained from the expression for $[X]$ by formally replacing $\bL$ by $q^2$.
\end{remark}

\smallskip

The expression of Proposition \ref{K0barM0n} for the Grothendieck classes 
$[\overline{M}_{0,n}]$ can also be interpreted as giving the counting of 
points over ``extensions" $\F_{1^m}$. In fact, the number of points over $\F_1$ 
can be obtained as the limit as $q\to 1$ of the function $N_X(q)$ 
that counts points over finite fields $\F_q$, possibly 
normalized by a power of $q-1$. 
The value $N_X(1)$ for a polynomially countable variety coincides 
with its Euler characteristic. 
Similarly, one can make sense of the number of points over
$F_{1^m}$ as the values $N_X(m+1)$, see Theorem 4.10 of \cite{CC2} 
and Theorem 1 of \cite{De}. 

\smallskip

\begin{proposition}\label{F1mbarM0n}
Let $p_{n,m}$ denote the number of points of 
$[\overline{M}_{0,n}]$ over $\F_{1^m}$. The generating function
$$ \varphi_m(t) =  \sum_{n\geq 1} p_{n,m} \frac{t^n}{n!} $$
is a solution of the differential equation
$$ (1+(m+1)\, t - (m+1)\, \varphi_m(t)) \, \varphi_m^\prime(t) = 1 + \varphi_m(t). $$
\end{proposition}

\proof Let $X$ be a smooth projective variety over $\Z$ whose class in 
the Grothendieck group can be written as $[X]=\sum_i a_i \bL^i$ 
with the $a_i$ non--negative integers. For all but finitely many primes
$p$ and $q=p^r$, the function that counts points of a finite fields $\F_q$
is then given by $N_X(q)=\sum_i a_i  q^i$. Thus, we obtain the values
$N_X(m+1)$ counting $F_{1^m}$-points by formally replacing $\bL$
with $m+1$ in the expression for the Grothendieck class.
\endproof

\medskip

\subsection{The moduli spaces $T_{d,n}$}

We consider here a family of varieties $T_{d,n}$ constructed in \cite{CGK},
which are natural generalizations of the moduli spaces $\overline{M}_{0,n}$.

\smallskip

We recall the construction of \cite{CGK} of $T_{d,n}$ as a 
family of varieties whose points parameterize stable $n$--pointed rooted trees of projective 
spaces $\P^d$. They generalize the moduli spaces $\overline{M}_{0,n}$, 
with the latter given by $T_{1,n}=\overline{M}_{0,n+1}$. These varieties are also closely
related to the Fulton--MacPherson compactifications $X[n]$ of configuration spaces, \cite{FM}
in the sense that, for any choice of a smooth complete variety $X$
of dimension $d$, one can realize $T_{d,n}$ in a natural way 
as a subscheme of $X[n]$.

\smallskip

\subsection{$n$--pointed rooted trees of projective spaces}\label{treePdSec}
A  graph $\tau$ consists of the data $(F_\tau, V_\tau,  \delta_\tau, j_\tau)$:
a set of flags (half-edges) $F_\tau$; a set of vertices $V_\tau$; 
 boundary map $\partial_\tau: F_\tau \to V_\tau$ 
that associates to each flag its boundary vertex;
and  finally the involution $j_\tau: F_\tau \to F_\tau$, $j_\tau^2=1$ that 
registers the matching of half--edges forming the edges of $\tau$. We consider here only graphs
whose geometric realizations are trees, 
i.e.~they are connected and simply connected.  
\\smallskip

A structure of rooted tree is defined by the choice of root tail $f_{\tau}\in F_{\tau}$,
$j(f_{\tau})=f_{\tau}$. Its vertex $v_{\tau}:=\delta (f_{\tau})$ also may be called the root.

\smallskip

We define the canonical orientation  on the rooted trees:  the root tail is oriented away 
from its vertex (so it is the {\it output});
all other flags are oriented towards the root vertex. The remaining tails are called inputs.

\smallskip

The output tail of a tree can be  grafted to an input tail of another tree.

\smallskip

We say that a vertex $v$ is a mother for a vertex $v'$ if $v'$ lies on an 
oriented path from $v$ to the root vertex $v_0$ and the oriented path from 
$v$ to $v'$ consists of a single edge. 

\smallskip

Given an oriented rooted tree $\tau$, we assign to each vertex $v\in V_\tau$
a variety $X_v\simeq \P^d$. To the unique outgoing tail  at $v$ we assign a choice
of an hyperplane $H_v \subset X_v$. To each incoming tail $f$ at $v$ we
assign a point $p_{v,f}$ in $X_v$ such that $p_{v,f}\neq p_{v,f'}$ for $f\neq f'$
and with $p_{v,f}\notin H_v$, for all  $f$ at $v$.

\smallskip

We think of an oriented rooted tree $\tau$, with $S_\tau$ the finite set of incoming 
tails of $\tau$ of cardinality $n$, as an $n$--ary
operation that starts with the varieties $X_{v_i}\simeq \P^d$ attached to the
input vertices $v_i$, $i=1,\ldots, m \leq n$, and glues the hyperplane 
$H_{v_i}\subset X_{v_i}$ to the exceptional divisor of the blowup of $X_{w_i}$ 
at the point $p_{w_i,f_i}$ where
$w_i$ is the target vertex of the unique outgoing edge of $v_i$ and $f_i$ is the
flag of this edge with $\partial(f_i)=w_i$, ingoing at $w_i$. The operation continues
in this way at the next step, by gluing the hyperplanes $H_{w_i}$ to the exceptional
divisor of the blowups of the projective spaces of the following vertex. At each vertex
that has an incoming tail, the corresponding variety acquires a marked point.
The variety obtained by this series of operation, when one reaches the root vertex, is the output
of $\tau$. It is endowed with $n$ marked points from the incoming tails and with 
a given hyperplane from the outgoing tail at the root.  In the terminology of 
\cite{CGK}, the output $X_\tau$ of an oriented rooted tree $\tau$ with $n$ incoming tails 
is a {\em $n$--pointed rooted tree of $d$--dimensional projective spaces}.

\smallskip

The {\em stability} condition for $X_\tau$ is the requirement that each component of
$X_\tau$ contains at least two distinct markings, which can be either marked points
or exceptional divisors. By Proposition 2.0.5 of \cite{CGK}, this condition 
is equivalent to the absence of nontrivial automorphisms of $\P^d$ fixing a hyperplane pointwise,
that is, translations and homotheties in $\A^d$.

\smallskip

Theorem 3.4.4 of \cite{CGK} defines the variety $T_{d,n}$ as the moduli space of 
$n$--pointed stable rooted trees of $d$--dimensional projective spaces. 

\medskip

\subsection{The class of $T_{d,n}$ and $\F_{1^m}$-points}\label{sec2.4}

In \cite{CGK}, the Poincar\'e polynomials of the varieties $T_{d,n}$ are computed,
generalizing the result of \cite{Man2} on the Poincar\'e
polynomial of the moduli spaces $\overline{M}_{0,n}$. 
Again, the classes $[T_{d,n}]$ in the Grothendieck ring can be computed with
the same technique, which shows that they satisfy the torification condition.
One also obtains the counting of points over $\F_{1^m}$.

\begin{proposition}\label{Prop2.4.1}
(1) For fixed $d$, the classes $[T_{d,n}]\in K_0(\cV_\Z)$ form a generating function
\begin{equation}\label{genTdnK0}
\psi(t) = \sum_{n\geq 1} [T_{d,n}] \frac{t^n}{n!} \, , 
\end{equation}
in $K_0(\cV_\Z)_\Q[[t]]$, which is the unique solution in $t+ t^2 K_0(\cV_\Z)_\Q[[t]]$
of the differential equation
\begin{equation}\label{eqdiffK0Tdn}
 (1 + \bL^d\, t - \bL \, [\P^{d-1}] \psi(t)) \, \psi^\prime(t) = 1 + \psi(t), 
\end{equation} 
where $[\P^{d-1}]  = \frac{\bL^d-1}{\bL -1}$. 

(2) The classes $[T_{d,n}]\in K_0(\cV_\Z)$ have a decomposition 
\eqref{K0torif} with non-negative coefficients.

(3) For a fixed $d$, denote by $p_{n,m}$ the number
of points of $T_{d,n}$ over $\F_{1^m}$ and form a generating function
$$ \eta_m(t) =  \sum_{n\geq 1} \frac{p_{n,m}}{n!} t^n.$$
This function  is a solution of the differential equation
$$ (1+(m+1)^d t - (m+1) \kappa_d(m+1) \eta_m) \eta_m^\prime = 1+\eta_m, $$
with $\kappa_d(q^2) = \frac{q^{2d}-1}{q^2 -1}$.
\end{proposition}

\smallskip

\proof (1) Theorem 5.0.2 and Corollary 5.0.3 of \cite{CGK} shows that,
for a fixed $d$ and for $n \geq 2$, the generating series
$$ \psi(q,t) = \sum_{n\geq 1} \frac{P_n(q)}{n!} \, t^n, $$
for the Poincar\'e polynomials $P_n(q):=P_{T_{d,n}}(q)$, with $P_1(q)=1$,
is the unique solution in $t+ t^2 \Q[q][[t]]$ to the differential equation
$$ (1 + q^{2d} t - q^2 \kappa_d(q^2) \psi) \partial_t \psi = 1 + \psi, $$
where $\kappa_d(q^2)$ is the Poincar\'e polynomial of $\P^{d-1}$. 
This result is obtained using the description of the varieties $T_{d,n}$
as iterated blowups, given in Theorem 3.6.2 of \cite{CGK}. 
The same construction of $T_{d,n}$, using the
blowup formula for the Grothendieck class, 
\begin{equation}\label{BlK0}
 [ {\rm Bl}_Y(X) ] = [X] + [Y] ([\P^{{\rm codim}_X(Y)-1}]-1), 
\end{equation} 
gives an analogous result for the classes. Namely, the relation
$$ P_{n+1}(q)= (\kappa_{d+1} + n q^2 \kappa_{d-1}) P_n(q)
+ q^2 \kappa_d \sum_{i+j=n+1, 2\leq i\leq n-1} \binom{n}{i} P_i(q) P_j(q) $$
satisfied by the Poincar\'e polynomials, as shown in \cite{CGK} is
replaced by the analogous relation for the Grothendieck classes
\begin{equation}\label{recK0Tdn}
[T_{d,n+1}] = ([\P^d] + n \bL\, [\P^{d-2}]) [T_{d,n}] + \bL \, [\P^d]\,
\sum_{i+j=n+1, 2\leq i\leq n-1} \binom{n}{i} [T_{d,i}]\, [T_{d,j}].
\end{equation}
These relations, for Poincar\'e polynomials and Grothendieck classes, 
respectively, can be seen from the inductive presentation of the Chow
group and the motive of $T_{d,n}$ given in \S 4 of \cite{CGK}, with
\eqref{recK0Tdn} following from the formula for the motive in Theorem 4.1.1 of
\cite{CGK}.

\smallskip

(2) The existence of a decomposition \eqref{K0torif} of $[T_{d,n}]$ 
with non-negative coefficients then follows from the fact that these
classes satisfy the recursive relation \eqref{recK0Tdn},
analogous to \eqref{recK0M0n}, which can be used to prove
the statement inductively, as in the case of $\overline{M}_{0,n}$.

\smallskip

(3) The counting $N_{T_{d,n}}(m+1)$ of $\F_{1^m}$-points
is obtained, as in the case of $\overline{M}_{0,n}$ 
by formally replacing $\bL$ with $m+1$ in the expression
for the Grothendieck classes or equivalently by replacing $q^2$ with
$m+1$ in the Poincar\'e polynomial.
\endproof

\bigskip

\section{Complemented subspaces and constructible sets}\label{ConstrSec}

In Borger's approach to $\F_1$-geometry via $\Lambda$-rings, \cite{Borg},
one has a notion of {\em complemented $\F_1$-points}. Namely, a sub-$\Lambda$-space 
$Y \subset X$ is {\em complemented} if the complement  $X\smallsetminus Y$
admits a $\Lambda$-space structure so that the map $X\smallsetminus Y \hookrightarrow X$
is a morphism of  $\Lambda$-spaces. In the case of toric varieties, with the
$\Lambda$-space structure determined by the torus orbits, the complemented
subspaces are unions of closures of torus orbits. In particular, the ``complemented $F_1$-points"
are the fixed points of the torus action, whose number equals the Euler characteristic.

\smallskip

The approach to $\F_1$-geometry via torifications is weaker than the approach
via $\Lambda$-rings. For example, as observed in Example 2.8 of \cite{Borg}, with
the exception of projective spaces, flag varieties do not admit a $\Lambda$-space
structure, though they certainly admit (non-affine) geometric torifications. However,
it is possible to consider an analogous notion of complemented subspaces in the
setting of torifications.

\medskip

\subsection{Complemented $\F_1$-points and torifications} 

Our use, in the previous section, of the decomposition into tori
of the class in the Grothendieck ring in order to count $\F_{1^m}$
points is based on thinking, as in \cite{CC2}, \cite{De}, of this
counting as being given by the values $N_X(m+1)$ of the
polynomial $N_X(q)$ counting points over $\F_q$. In terms of
Grothendieck classes, we obtained the counting of $\F_{1^m}$-points
as
\begin{equation}\label{countF1m}
\# X(\F_{1^m})  = [ X ] |_{\bT = m} = \sum_k a_k m^k, 
\end{equation}
by formally replacing the variable $\bT$ with $m$ in the expression
$[X]=\sum_k a_k \bT^k$, with $a_k\geq 0$ for the Grothendieck
class. The case of points over $\F_1$ corresponds to $m=0$, with
$\# X(\F_1)  = [ X ]|_{\bT =0} = a_0=\chi(X)$. 
Essentially, this means that, for a variety $X$ with a torification,
only the zero dimensional points contribute $\F_1$-points,
while each $k$-dimensional torus $T^k$ of the torification with $k>0$ 
contributes $m$ points over $\F_{1^m}$ for each $m\geq 1$. This is 
related to the general philosophy that the extensions $\F_{1^m}$
are related to actions of the groups $\mu_m$ of $m$-th roots of unity, 
see \cite{KS} and more recently \cite{CCM}, \cite{Man1}.

\smallskip

This leads to a natural generalization of the notion of complemented
$\F_1$-points in the context of torifications. The counting formula 
\eqref{countF1m} implies that, according to this notion, $\F_1$-points
are points of $X$ such that the Grothendieck class of the complement
of these points still admits a decomposition into tori with non-negative
coefficients. At the level of geometric torifications, it is natural therefore
to introduce a stronger notion of complemented points as follows.

\begin{definition}\label{F1compl}
Let $X$ be a variety over $\Z$, with a geometric torification.
A finite set of points $S$ is (strongly, ordinarily, weakly) {\em complemented} 
if the complement $X\smallsetminus S$ also has a 
geometric torification, such that the inclusion $X\smallsetminus S
\hookrightarrow X$ is a (strongly, ordinarily, weakly) torified morphism.
\end{definition}

For example, if we consider $\P^1$ with a torification given by the choice
of two points, each of these two points is a strongly complemented $\F_1$-point,
while any other point would be an ordinarily complemented $\F_1$-point, since
the complement can be torified and the inclusion becomes a torified morphism
after composing with an isomorphism of $\P^1$. 

\medskip

\subsection{Complemented torifications}

Similarly, one has a notion of {\em complemented subspace} in a 
torified variety. Torifications behave well with respect to blowups
along complemented subspaces. 

\begin{definition}\label{F1complY}
Let $X$ be a variety over $\Z$ with a geometric torification. 
A subvariety $Y \subset X$ is (strongly, ordinarily, weakly) {\em complemented} 
if both $Y$ and the complement $X\smallsetminus Y$ have a 
geometric torifications, so that the inclusions $Y\hookrightarrow X$ and 
$X\smallsetminus Y \hookrightarrow X$ are (strongly, ordinarily, weakly) 
torified morphisms.
\end{definition}

On a variety that has a geometric torification compatible with
a $\Lambda$-structure, the complemented condition for  
sub-$\Lambda$-spaces of \cite{Borg} implies
the strong form of complementation of Definition \ref{F1complY}.

\begin{example}\label{diagex}
Consider $\P^n$ with a torification $\cT=\{ T_i \}$ and $\P^n \times \P^n$ 
with the torification $\{ T_i \times T_j \}$. The diagonal $\Delta \subset \P^n \times \P^n$
is weakly complemented but not ordinarily nor strongly complemented. In the big
cell $\A^n \times \A^n$ with the product torification, the diagonal is ordinarily
complemented but not strongly complemented.
\end{example}

We have the following behavior of torifications with respect to blowups.

\begin{proposition}\label{complblowup}
Let $X$ be a variety over $\Z$ with a geometric torification and let $Y\subset X$
be a (strongly, ordinarily, weakly) {\em complemented} subspace. Then the
torifications of $Y$ and of $X\smallsetminus Y$ for which the inclusions
are (strongly, ordinarily, weakly) torified morphisms determine a geometric
torification of the blowup ${\rm Bl}_Y(X)$ of $X$ along $Y$, such that the morphism
$\pi: {\rm Bl}_Y(X) \to X$ is (strongly, ordinarily, weakly) torified.
\end{proposition}

\proof It suffices to show that the (strongly, ordinarily, weakly) 
compatible torifications of $X$, $Y$ determine a torification of
the exceptional divisor of the blowup, since the complement
is then torified by the torification of $X\smallsetminus Y$. Thus,
we consider the projectivized normal bundle $\P(N_X(Y))$.
The restriction of the bundle over the tori of the torification of
$Y$ is trivial, hence $\P(N_X(Y))$ can be torified by the products 
of the tori in the torification of $Y$ with the tori in a torification of
$\P^{{\rm codim_X(Y)}-1}$.
\endproof

The blowup operation does not behave well with respect to
geometric torifications in the non-complemented case. For example,
the blowup of a 2-dimensional torus at a point does not have a 
torification compatible with the blowup morphism, even in the
weak sense.

\medskip

\subsection{Geometric torifications of stable curves of genus zero}

The fibers of the forgetful morphism $\overline{M}_{0,n+1} \to \overline{M}_{0,n}$
consist of the stable curves of genus zero. We show that these admit geometric 
torifications. In general, these torifications are neither regular nor affine.

\begin{lemma}\label{Ctor}
Let $C$ be a stable (pointed) curve of genus zero. A choice of a rooted tree,
for a torification of $\P^1$, and of a point in each other component
determines a geometric torification of $C$. 
\end{lemma}

\proof A geometric torification of $C$
is obtained by identifying $C$ with a tree of $\P^1$'s with marked points,
choosing a root vertex and a torification of the $\P^1$ at the root, for which
the two zero dimensional tori are away from the intersection points with
other components (and from the additional marked points), and then at each 
adjacent vertex choosing a torification given by a torification of the 
complement of the intersection point with the $\P^1$ of the root vertex, 
with the zero dimensional torus chosen away from the intersection points 
with other $\P^1$'s (and away from the additional marked points), and so on.
\endproof

These torifications reflect the decomposition into tori of the
class $[C]$ in the Grothendieck ring: a tree of
$\P^1$'s with $N$ vertices has class $[C]=N \bT+ N+1$.

\medskip

\subsection{Marked points, $\F_{1^m}$-points, and constructible sets}

When considering stable curves of genus zero with $n$ marked points,
the geometric torifications considered in Lemma \ref{Ctor} do not
reflect the presence of the marked points, as these are not part of the torification. 
However, in order to descend to $\F_1$ the notion of marked points, we
need to ensure that the stable curve $C$ has enough points, possibly
after passing to an extension to some $\F_{1^m}$.

\begin{lemma}\label{CptsF1m}
If $C$ is a tree of $\P^1$'s with $N$ vertices, then the number of points
of $C$ over $\F_{1^m}$ is $N(m+1)+1$. 
\end{lemma}

\proof This follows from the same argument used before, computing
the number of points over $\F_{1^m}$ from the Grothendieck class
by formally replacing $\bL=\bT+1$ with $m+1$. 
\endproof

In particular, in the case of a single $\P^1$ the number of
points over $\F_{1^m}$ is $m+1$. Thus, when we consider 
stable curves with genus zero with $n$ marked
points, one should work with an extension $\F_{1^m}$ with
$m\geq n-1$. Passing to an extension in this way is 
necessary in order to have morphism induced by the
action of $S_n$ that permutes points. We will discuss
more precisely the nature of such morphisms below.
This phenomenon is similar to what happens in \cite{CC1},
where Chevalley groups define varieties over
the extension $\F_{1^2}$.

\smallskip

Notice that, even after considering $\F_{1^m}$-points so as
to ensure the existence of the correct number of marked points, 
one needs to work with points that are not necessarily
complemented. 

\smallskip

In the usual case of algebraic varieties, the complement
of an algebraic variety inside another may not determine
an algebraic variety, but a {\em constructible set}. When
we consider stable curves of genus zero with marked points,
the fact that the points are not complemented means that 
the complement does not define an $\F_1$-variety. However,
one expects that it will define an $\F_1$-constructible set,
in an appropriate sense. We show in the rest of this
section how one can define a suitable notion of 
$\F_1$-constructible sets, by relaxing the notion of
geometric torification, while retaining intact the decomposition
of the class in the Grothendieck ring.

\bigskip

\section{Constructible sets over $\F_1$ and torifications}\label{ConstrSec2}

The need to consider points that are not complemented in the
case of the torifications of stable curves of genus zero, suggests
that one should take into consideration a further level of
structure that lives in between the coarse condition about the
decomposition of the Grothendieck class into a sum of tori with
non-negative coefficients and the geometric torifications, which
allows for a larger class of complemented subspaces and provides
a suitable notion of $\F_1$-constructible set.

\smallskip

Starting from the observation that torifications behave well with
respect to products and disjoint unions, but not with respect to
complements, one can introduce a weaker notion of {\em constructible 
torification}, which is similar conceptually to the usual way of passing 
from a semigroup to a group. Recall that every constructible subset $C$
of an algebraic variety determines a class in the Grothendieck ring
$K_0(\cV)$ of varieties.

\begin{definition}\label{CF1sets}
Let $\cC_{\F_1}$ be the class of constructible sets over $\Z$ that
can be obtained, starting from $\bG_m$ through the operations 
of products, disjoint unions, and complements. Let $X$ be a 
constructible set over $\Z$. A constructible torification
of $X$ is a morphism of constructible sets $e_X : C \to X$ from an element $C\in \cC_{\F_1}$ 
to $X$ such that the restriction of $e_X$  to each component of $C$ is an immersion 
and $e_X$ induces a bijections of $k$--points, $e_X(k): T(k) \to X(k)$, for every field $k$.
An $\F_1$-constructible set is a constructible set over $\Z$ together with a
constructible torification, with the property that the class $[X]$ in the Grothendieck
ring of varieties $K_0(\cV_\Z)$ has a decomposition \eqref{K0torif} in classes of tori,
with non-negative coefficients.
\end{definition}

The class $\cC_{\F_1}$ considered above includes all tori, as well as complements
of disjoint unions of tori inside other tori, products of such sets, and so on. 
$\F_1$-constructible sets are built out of these building blocks, with the 
requirement that the positivity condition on the tori decomposition of the
Grothendieck class holds.

As in the case of geometric torifications, one assigns an equivalence relation
between constructible torifications that corresponds to defining the same
structure of $\F_1$-constructible set. This can be done in a strong, ordinary
and weak form, following the analogous definitions for geometric torifications.

\begin{definition}\label{parTmor}
Let $X$ and $Y$ be constructible sets over $\Z$, endowed with 
constructive torifications. A morphism $f: X \to Y$ is a strong 
morphism of $\F_1$-constructible sets if for each component 
$C_j$ of the constructive torification of $Y$, $f^{-1}(C_j)$ is a disjoint
union of components of the partial torification of $X$. 
Let $X$ be a variety over $\Z$. Two constructive torifications of 
$X$ are strongly equivalent if the identity on $X$ is a strong 
morphism of $\F_1$-constructible sets. They are ordinarily equivalent
if there is an isomorphism of $X$ that is a strong 
morphism of $\F_1$-constructible sets and they are weakly equivalent 
if there are decompositions $\{ Z_k \}$ and $\{ Z'_k \}$ of $X$ compatible 
with the constructive torifications, and isomorphisms
$\phi_k: Z_k \to Z'_k$ that are strong 
morphisms of $\F_1$-constructible sets.
An ordinary morphism of $\F_1$-constructible sets is a morphism $f: X \to Y$
such that $\psi\circ f \circ \phi$ is a strong morphism of $\F_1$-constructible sets,
for some $\phi$ and $\psi$ isomorphisms of $X$ and $Y$.
A weak morphism of $\F_1$-constructible sets is a morphism $f: X \to Y$ such that
$\psi_k \circ f \circ \phi_\ell$ is a strong morphism, where $\psi_k$ and
$\phi_\ell$ are isomorphisms of pieces of decompositions of $Y$ and
$X$, respectively, compatible with the constructible torification.
\end{definition}

The following result on blowups for constrictive torifications will
be useful later.

\begin{lemma}\label{constrTblowup}
Let $X$ be a variety over $\Z$ with a constructible torification.
Let $Y\subset X$ be a closed subvariety, such that $X\smallsetminus Y$ has a 
constructible torification and $Y$ has a geometric torification and the
inclusions are (strong, ordinary, weak) morphisms of constructibly
torified spaces. Then these torifications determine a constructible 
torification of the blowup ${\rm Bl}_Y(X)$ so that the map
$\pi: {\rm Bl}_Y(X) \to X$ is a (strong, ordinary, weak) morphism 
of constructibly torified spaces.
\end{lemma}

\proof The argument is as in Proposition \ref{complblowup}.
The exceptional divisor, which we identify with $\P(N_X(Y))$
has a geometric torification, since it is trivial when restricted
to the tori of the geometric torification of $Y$. 
The constructible torification of $X\smallsetminus Y$
extends the torification of the exceptional divisor to a
constructible torification of ${\rm Bl}_Y(X)$.
\endproof

For this construction to extend
to the case where the blowup locus $Y$ has a constructible
torification, one would need to ensure that the bundle $\P(N_X(Y))$
is trivial when restricted to the components $C_i\in \cC_{\F_1}$
of the decomposition of $Y$. This is the case, for instance, when
the complements of unions of tori inside other tori in the
sets $C_i$ extend to actual (not necessarily disjoint) tori in $Y$.

\medskip

\subsection{Categories of constructible torifications}

As in the case of geometric torifications, 
the different notions of morphisms of constructible torifications give
rise to different categories.

\begin{proposition}\label{catsCTor}
There are categories $\cC\cT^s \subset \cC\cT^o \subset \cC\cT^w$ where
the objects, ${\rm Obj}(\cC\cT^s)={\rm Obj}(\cC\cT^o)={\rm Obj}(\cC\cT^w)$,
are pairs $(X_\Z, \cC)$, with $X_\Z$ a constructible set over $\Z$ and $\cC=\{ C_i \}$
is a constructible torification of $X_\Z$, in the sense of Definition \ref{CF1sets}.
Morphisms in $\cC\cT^s$ are strong morphisms of
constructibly torified spaces; morphisms in $\cG\cT^o$ are ordinary
morphisms of constructibly torified spaces; 
morphisms in $\cG\cT^w$ are weak 
morphisms of constructibly torified spaces.
\end{proposition}

\proof Strong morphisms of constructibly torified spaces are
as in Definition \ref{parTmor}. 
The torified condition is preserved by composition.
Morphisms in $\cC\cT^o$ are arbitrary compositions of  
strong morphisms and isomorphisms of $\Z$-constructible sets,
hence composition is also well defined. These are the
ordinary morphisms of constructibly torified spaces, as
in Definition \ref{parTmor}.
Morphisms in $\cC\cT^s$ are also morphisms
in $\cC\cT^o$, but in general not conversely.
Similarly, morphisms in $\cC\cT^s$ are arbitrary compositions
of strong morphisms and weak equivalences, in the sense
of Definition \ref{parTmor}. Composition
is well defined and morphisms in $\cC\cT^s$ and in $\cC\cT^o$ 
between any pair of objects are proper subsets of morphisms
in $\cC\cT^w$ between the same objects.
\endproof

\bigskip

\section{Constructible torifications of moduli spaces}\label{ctorSec}

We apply the notion of $\F_1$-constructible sets introduced
above to define $\F_1$-structures on the moduli spaces
$\overline{M}_{0,n}$ and on their generalizations $T_{d,n}$.

\subsection{Constructible torification of $\overline{M}_{0,n}$}

As we have seen, stable curves of genus zero with marked points are
$\F_1$-constructible sets. It is therefore natural to seek a realization of
the moduli spaces $\overline{M}_{0,n}$ in $\F_1$-geometry also as
$\F_1$-constructible sets. We show that the moduli spaces 
$\overline{M}_{0,n}$ have a constructible torification, underlying the decomposition
of the Grothendieck class into tori. 

\smallskip

\begin{theorem}\label{ConstrTbarM0n}
The moduli spaces $\overline{M}_{0,n}$ are $\F_1$-constructible
sets with a constructible torification determined by the choice of
a constructible torification of $\P^1$ minus three points.
\end{theorem}

\proof For $n\geq 4$, we identify $M_{0,n}$ with the complement of the diagonals
in the product of $n-3$ copies of $\P^1$ minus three points. The complement of
three points in $\P^1$ is an $\F_1$-constructible set, with a constructible torification 
given by two points and the complement of one point in a 1-dimensional torus. The
product of $n-3$ copies of $\P^1 \smallsetminus \{ 0,1,\infty \}$, in turn has the product 
constructible torification. When we remove the diagonals, this corresponds to taking
complements of sets in the class $\cC_{\F_1}$ inside other sets in the same class,
hence we still obtain a set in $\cC_{\F_1}$. This does {\em not} define a structure
of $\F_1$-constructible set  on the open stratum $M_{0,n}$ by itself, because the positivity
condition on the class $[ M_{0,n} ]$ is not satisfied. However, we consider the
constructible torification of $M_{0,n}$ together with the constructible torifications
obtained in this same way on all the products $\prod_i M_{0,k_i+1}$ with 
$\sum_i k_i=n$, of the lower dimensional strata and we obtain a 
constructible torification of $\overline{M}_{0,n}$, which also satisfies the
positivity condition on the Grothendieck class, by Proposition \ref{K0barM0n}.
Thus, $\overline{M}_{0,n}$ is an $\F_1$-constructible set.
\endproof

\medskip

\subsection{Constructible torification of $T_{d,n}$}\label{Sec2}

We extend here the construction of geometric torifications described above 
for $\overline{M}_{0,n}$ to the case of the $T_{d,n}$ of \cite{CGK}.

\smallskip

We have see that $n$-pointed stable curves of genus zero define
$\F_1$-constructible sets by a choice of a geometric torification of the
underlying tree of $\P^1$'s and the constructible torification given by
taking the complement of the marked points in the geometric torification.
One has an analogous construction for the $n$-marked stable trees 
of projective spaces described in \S \ref{treePdSec}. 

\begin{proposition}\label{treePdtor}
Let $\Gamma$ be an $n$-marked stable tree of projective spaces.
Then $\Gamma$ defines an $\F_1$-constructible set with a constructible
torification determined by the choice of a torification of $\cT_{\P^d}$.
\end{proposition}

\proof Given an oriented rooted tree $\tau$, at the root vertex $v_0$
we assign a choice of a torification $\cT_{\P^d}=\{ T_i \}$ of $X_{v_0}\simeq \P^d$,
with a compatible torification of the hyperplane $H_{v_0}$ at the unique
outgoing tail at $v_0$. We then replace the torification $\cT_{\P^d}$ of $X_{v_0}$
by a constructible torification of the complement in $X_{v_0}$ of the marked
points corresponding to all the incoming tails at $v_0$, by replacing the tori $T_i$ of
the torification that contain a subset $\{ p_{v_0, f, i} \}$
of the marked points with the sets $C_i \in \cC_{\F_1}$
given by the complements $C_i=T_i\smallsetminus \cup p_{v_0, f, i}$. By construction,
the marked points on $X_{v_0}$ are distinct and not contained in the 
hyperplane $H_{v_0}$. We then consider the blowup of $X_{v_0}$ at the marked
points corresponding to incoming edges of the tree. 
By Lemma \ref{constrTblowup}, the constructible torification of the complement
of the marked points determines a constructible torification on the blowup. One
then considers the adjacent vertices $v_i$, and glues the hyperplane $H_{v_i}$
to the exceptional divisor of the blowup of $X_{v_0}$ at the tail mark of the edge 
from $v_i$ to $v_0$, so as to match the torification of $H_{v_i}$ with a piece of
the torification of the exceptional divisor. One continues in this way for all the
other vertices. 
\endproof

We also obtain a constructible torification of $T_{d,n}$, by a construction
similar to what we have for $\overline{M}_{0,n}$. 

\begin{theorem}\label{TdnConstrT}
The moduli spaces $T_{d,n}$ are $\F_1$-constructible
sets with a constructible torification determined by the choice of
a constructible torification of $\A^d$ minus two points.
\end{theorem}

\proof The open stratum $TH_{d,n}$ of $T_{d,n}$ is the configuration space
of $n$ distinct points in $\A^d$ up to translation and homothety, or equivalently
of all embeddings of a hyperplane $H$ and $n$ distinct points not on the
hyperplane in $\P^d$, up to projective automorphisms that pointwise fix $H$.
Fixing two points in $\A^d$ suffices to fix the symmetries, since fixing the origin 
${\bf 0}=(0,\ldots,0)$ eliminates translations and fixing another point, for instance
${\bf 1}=(1,\ldots, 1)$, takes care of homotheties. 
Thus, we can identify the open stratum with
\begin{equation}\label{THdn}
TH_{d,n} \simeq (\A^d \smallsetminus \{ {\bf 0}, {\bf 1} \})^{n-2} \smallsetminus \Delta, 
\end{equation}
the complement of all the diagonals $\Delta$ in the product of $n-2$ copies
of $\A^d$ minus two points. In the case where $d=1$, this gives back the usual
description of $M_{0,n+1}=TH_{d,n}$ as the complement of the diagonals
in a product of copies of $\P^1\smallsetminus \{ 0,1,\infty \}$. A choice
of a constructible torification on the complement of two points in $\A^d$ then
determines constructible torifications on the products, on the diagonals,
and on the complements. As shown in \cite{CGK}, the compactification $T_{d,n}$
has boundary components isomorphic to products $T_{d,n_1}\times T_{d,n_2}$
with $n_1 + n_2=n$. Thus, by considering constructible torifications on all the
open strata $\prod_i TH_{d,n_i}$ with $\sum_i n_i =n$, we obtain a constructible
torification of $T_{d,n}$. The condition of positivity of the Grothendieck class
is not in general satisfied by the individual $TH_{d,k}$ and their products, but
it is satisfied by $T_{d,n}$ itself, because of Proposition \ref{Prop2.4.1}. Thus,
the moduli spaces $T_{d,n}$ have a structure of $\F_1$-constructible sets.
\endproof

\bigskip

\section{Morphisms and operad structure}\label{opSec}

We show that the constructions of torifications described in the previous section
are compatible with the operad structures.

\smallskip

Let an operad ${\mathcal P}$, in the symmetric monoidal category
of varieties over $\Z$ with Cartesian product be given. Its descent data to $\F_1$
consist of affine torifications such that the 
composition operations
$$ {\mathcal P}(n) \times {\mathcal P}(m_1) \times \cdots \times {\mathcal P}(m_n)
\to {\mathcal P}(m_1+\cdots + m_n) $$
and the structure actions of symmetric groups
are morphisms of affinely torified varieties.
This is a ``tourist class''  description of \cite{Man3}. A more systematic treatment requires
the explicit introduction of a category of labeled graphs as in \cite{BoMan}.

\smallskip

\subsection{Categories of trees and operads}\label{sec3.1}

We consider a category $\Gamma$, whose
objects are finite disjoint unions of oriented rooted trees.
The morphisms are generated by edge contractions and
graftings.  The grafting of an oriented tree $\tau$ to another oriented tree $\sigma$
is realized by the morphism $h: \tau \amalg \sigma \to \tau \#_{v_0,w} \sigma$,
where the involution $j_h$ matches the outgoing tail of the root vertex $v_0$ 
of $\tau$ with an ingoing tail of a vertex $w$ of $\sigma$. The edge contractions
are given by morphisms $h_e: \tau \to \tau/e$, where the edge $e$ is a $j_\tau$-orbit
$e=\{ f, f'\}$ of flags $f,f'\in F_\tau$ and $F_\tau \setminus h_e^F(F_{\tau/e})= \{ f,f' \}$,
$h_{e,F}^{-1}: h_e^F(F_{\tau/e}) \to F_{\tau/e}$ the identity and $h_{e,V}: V_\tau \to V_{\tau/e}$
mapping $\partial_\tau(f)$ and $\partial_\tau(f')$ to the same vertex in $\tau/e$.

\smallskip

It is shown in Section IV.2 of \cite{Man3} that the datum of
an operad is equivalent to a monoidal functor ${\mathcal M}$ 
from a category of trees (forests) with the symmetric monoidal structure given 
by disjoint union and morphisms generated by graftings and 
edge contractions, to a symmetric monoidal category 
$({\mathcal C},\otimes)$, with the condition that 
$$ {\mathcal M}(\tau) = \otimes_{v\in V_\tau} {\mathcal M}(\tau_v), $$
where $\tau_v$ is the star of the vertex $v$, see Proposition IV.2.4.1 of \cite{Man3}.
The operad composition is identified with the image ${\mathcal M}(\psi)$ of the morphism
$\psi$ that assigns to a disjoint union of corollas $\tau \amalg \tau_1 \cdots \amalg \tau_n$ 
the corolla obtained by first grafting the outgoing tails 
of the component $\tau_k$ to the $k$-th ingoing tail of $\tau$ and then contracting all
the edges.

\subsection{Operad morphisms of $\overline{M}_{0,n}$}

We now consider the composition maps that give the operad structure of
the moduli spaces $\overline{M}_{0,n}$ and see that these are also compatible
with the structure of $\F_1$-constructible sets described above.

\begin{theorem}\label{thm3.3.1}
Let ${\mathcal M}(n)=\overline{M}_{0,n+1}$.
The composition morphisms of the operad 
$$ {\mathcal M}(n) \times {\mathcal M}(m_1) \times \cdots \times {\mathcal M}(m_n) \to 
{\mathcal M}(m_1+\cdots +m_n) $$
are strong morphisms of constructibly torified spaces, with respect to 
the constructible torifications of Theorem \ref{ConstrTbarM0n}. Thus, the
operad ${\mathcal M}(n)$ descends to an operad of $\F_1$-constructible sets,
in the category $\cC\cT^s$ of Proposition \ref{catsCTor}.
\end{theorem}

\proof The constructible torification of $\overline{M}_{0,n}$, obtained as in
Theorem \ref{ConstrTbarM0n}, is built out of constructible torifications of
the open strata $M_{0,k}$ and their products, so that one has a family of
compatible constructible torifications on the open strata 
$\prod_k M_{0,n_k+1}$ with $\sum_k n_k=n$. This implies that the
inclusions of the boundary strata $\prod_k \overline{M}_{0,n_k+1}$
are compatible with the constructible torifications, hence they are 
strong morphisms of $\F_1$-constructible sets.
\endproof

The symmetric group $S_n$ acts on $\overline{M}_{0,n}$ by
permuting the marked points. 

\begin{proposition}\label{Snact}
The elements of the symmetric group $S_n$ act on 
$\overline{M}_{0,n}$ as ordinary 
morphisms of $\F_1$-constructible sets, that is, morphisms
in the category $\cC\cT^o$ of Proposition \ref{catsCTor}.
\end{proposition}

\proof
The constructible torification of $\overline{M}_{0,n}$ described in
Theorem \ref{ConstrTbarM0n} is obtained from a constructible
torification of $\P^1\smallsetminus \{ 0, 1, \infty \}$ given by the
points $0$, $\infty$ as zero-dimensional tori and a the complement
$\bG_m \smallsetminus \{ 1 \}$ of a zero-dimensional torus in a 
1-dimensional torus as the remaining piece of the decomposition
into sets of $\cC_{\F_1}$. The action of an element $\sigma\in S_n$ 
on $M_{0,n}$ is a  permutation of the $n$ marked points and is
therefore given by an isomorphism of
$M_{0,n}$ that sends this choice of a constructible torification into a 
different choice, obtained by a different initial choice of constructible
torification of $\P^1\smallsetminus \{ \sigma(0), \sigma(1), \sigma(\infty) \}$
and a permutation of the factors in the poduct 
$(\P^1\smallsetminus \{ 0, 1, \infty \})^{n-3}$. Thus, the permutation
group $S_n$ acts on $\overline{M}_{0,n}$ by ordinary (not strong) 
morphisms of $\F_1$-constructible sets.
\endproof

\medskip

\subsection{Operad morphisms of $T_{d,n}$}

The varieties $T_{d,n}$ have natural morphisms defining an
operad structure that generalizes the operad of $\overline{M}_{0,n}$.
We use here, for convenience, the notation $T_{d,S}$, with $S$ 
the set of marked points, with $\# S=n$.

\medskip

\begin{theorem}\label{thm3.2.1}
For each fixed $d\geq 1$, there are morphisms  
of the following form, which determine an operad ${\mathcal T}_d$:

(1) isomorphisms: $T_{d,S} \stackrel{\simeq}{\to} T_{d,S'}$ for $S'\stackrel{\simeq}{\to} S$,
functorial wrt the bijections of labelling sets.

(2) embeddings: $T_{d,S'} \times T_{d,S\setminus S' \cup \{ \star \}} \hookrightarrow T_{d,
S}$, for $S'\subset S$ with $\# S' \geq 2$.

(3) forgetful morphisms: $T_{d,S} \to T_{d,S'}$ for $S'\subset S$ with
$\# S' \geq 2$.

These morphisms satisfy the standard identities.
\end{theorem}

\smallskip

\proof
The existence of morphisms of the form (1) is clear by construction. 
The cases (2) and (3) follow from the boundary stratification of these varieties 
constructed in Theorem 3.3.1 of \cite{CGK}. In fact, the boundary of a variety 
$T _{d,S}$ is given by smooth normal crossings divisors: given any proper
subset $S'\subset S$, there is a nonsingular divisor $T_{d,S}(S')\subset T_{d,S}$.
These divisors meet transversely and the only non-empty intersections 
$T_{d,S}(S_1)\cap \cdots \cap T_{d,S}(S_r)$ happen when the sets $S_k$ are
nested (each pair is either disjoint or one is a subset of the other). The
divisors satisfy $T_{d,S}(S')\simeq T_{d,S'} \times T_{d,S\setminus S' \cup \{ \star \}}$.
This gives the morphisms (2) coming from the inclusion of the strata. In terms of
morphisms of oriented rooted trees, these correspond to the morphisms 
that graft the outgoing tail of the first tree with set of incoming tails identified with $S'$
to the incoming tail marked by $\star$ in the second tree. The forgetful
morphisms (3) come from the construction of $T_{d,S\cup \{ \star \}}$ from
$T_{d,S}$ via a sequence of iterated blowups, as in Theorem 3.3.1 of \cite{CGK}.
The composition of the projections of this sequence of blowups give the
forgetful morphism $T_{d,S\cup \{ \star \}} \to T_{d,S}$.  In terms of rooted
trees of projective spaces, these correspond to forgetting some of the marked
points and contracting the resulting unstable components.

\smallskip

Using the functorial characterization of operads given in
Proposition IV.2.4.1 of \cite{Man3}, let $({\mathcal C},\otimes)$ 
be the symmetric monoidal category of algebraic varieties with 
the Cartesian product, and let $(\Gamma, \amalg)$ be the category 
of oriented rooted forests with disjoint union. The embeddings
of the strata determine the morphisms ${\mathcal T}_d(\psi)$,
where $\psi$ is the morphism of oriented rooted trees that
assigns to a disjoint union of oriented corollas $\tau \amalg \tau_1 \cdots \amalg \tau_n$,
where each corolla has only one outgoing tail, the corolla obtained by first grafting 
the outgoing tails of the component $\tau_k$ to the $k$--th ingoing tail of $\tau$ 
and then contracting all the edges. This assignment determines the operad 
composition operations
$T_{d,S} \times T_{d,S_1} \times \cdots \times T_{d,S_n} \to T_{d, S_1\cup \cdots \cup S_n}$,
where $n=\# S$ is the number of incoming tails of the trees of projective spaces
parameterized by $T_{d,S}$. 
\endproof

\smallskip

\begin{theorem}\label{thm3.2.2}
The operad composition operations
\begin{equation}\label{composeTdn}
T_{d,k} \times T_{d,n_1} \times \cdots \times T_{d,n_k} \to T_{d, n_1+ \cdots + n_k}
\end{equation}
are strong morphisms of constructible torifications, hence they define strong
morphisms of $\F_1$-constructible sets, that is, morphisms in the category
$\cC\cT^s$ of Proposition \ref{catsCTor}. 
\end{theorem}

\smallskip

\proof  We consider the structure of $\F_1$-constructible sets on the moduli spaces $T_{d,n}$ given
by the constructible torification obtained as in Theorem \ref{TdnConstrT}. 
Since this is built as a collection of compatible constructible torifications on
all the boundary strata of $T_{d,n}$ we see that the operad
composition operations \eqref{composeTdn}, which are 
obtained from the morphisms of type (2) of Theorem \ref{thm3.2.1}, are inclusions
of boundary strata, hence compatible with the constructible torification.
\endproof

The result above accounts for the morphisms of type (2) in Theorem \ref{thm3.2.1}.
The morphisms of type (1) and (3) also determine morphisms of constructible torifications

\smallskip

\begin{proposition}\label{13morphTdn}
Morphisms of type (1) in
Theorem \ref{thm3.2.1} are ordinary (not strong) morphisms of $\F_1$-constructive sets.
Morphisms of type (3) of Theorem \ref{thm3.2.1}  are strong morphisms
of $\F_1$-constructible sets, that is, morphisms in the category
$\cC\cT^s$ of Proposition \ref{catsCTor}. 
\end{proposition}

\proof The case of morphisms of type (1) is analogous to the case of
the morphisms permuting the marked points of $\overline{M}_{0,n}$,
and for the same reason they are ordinary (not strong) $\F_1$-morphisms.
Morphisms of type (3) are the forgetful morphisms $T_{d,n+1} \to T_{d,n}$
that forget one of the marked points and contract the unstable components. 
The restrictions of these projection maps to the open strata 
$TH_{d,n}$ and $\prod_i TH_{d,n_i}$, with $\sum_i n_i=n$ and $n_i \geq 2$,
are given by projections $(\A^d\smallsetminus \{ {\bf 0}, {\bf 1} \})^{n_j-2}\smallsetminus \Delta
\to  (\A^d\smallsetminus \{ {\bf 0}, {\bf 1} \})^{n_j-3}$ on one of the factors with $n_j>2$.
By construction of the constructible torifications on the $T_{d,n}$ given in Theorem
\ref{TdnConstrT}, these projections are morphisms of $\F_1$-constructible sets.
\endproof

\bigskip

\section{Moduli spaces and wonderful compactifications}\label{wondSec}

Another approach to defining  $\F_1$-structures on the
operads of the moduli spaces $\overline{M}_{0,n}$ and
of the $T_{d,n}$ is based on the construction of the moduli spaces
$T_{d,n}$ as iterated blowups and their relation to the
Fulton--MacPherson compactifications. 

\medskip
\subsection{Moduli spaces $T_{d,n}$ and Fulton--MacPherson spaces}

Let us denote by $X[S]$ the Fulton--MacPherson space. We 
describe its construction in terms of iterated blowups, following the
general construction for graph configuration spaces used in 
\cite{Li1}, \cite{Li2} and in \cite{CeMar}, in the special case of the 
complete graph. One starts with the product $X^S$ of $n=\# S$ copies
of $X$ and considers all diagonals $\Delta_{S'}\subset X^S$ for
all subsets $S'\subseteq S$, given by $\Delta_{S'}=\{ x\in X^S \,|\,
x_i =x_j, \,\, \forall  i,j \in S' \}$. Upon identifying the subset $S'\subset S$
with the set of vertices of a subgraph $\Gamma_{S'} \subset \Gamma_S$,
where $\Gamma_S$ is the complete graph on the set $S$ of vertices,
the diagonal $\Delta_{S'}$ is identified with a product $X^{V_{\Gamma_S/\Gamma_{S'}}}$
where the quotient graph $\Gamma_S/\Gamma_{S'}$ is obtained by identifying all of
$\Gamma_{S'}$ to a single vertex. Consider the set ${\mathcal G}_S$ 
of all subgraphs $\Gamma_{S'}$ that are biconnected (that is, that cannot be 
disconnected by removing the star of any one vertex) and choose an ordering
${\mathcal G}_S = \{ \Gamma_{S'_1}, \ldots, \Gamma_{S'_N} \}$ such that
if $S_i \supseteq S_j$ then the indices are ordered by $i\leq j$. 
By dominant transform of a subvariety under a blowup one means  
the proper transform if the variety is not contained in the blowup locus
and the inverse image otherwise (Definition 2.7 of \cite{Li1}).
It was shown in Theorem 1.3 and Proposition 2.13 of \cite{Li1} (see also
Proposition 2 of \cite{CeMar}) that the sequence of blowups $Y^{(k)}$
with $Y^{(0)}=X^S$ and $Y^{(k+1)}$ obtained by blowing up $Y^{(k)}$
along the dominant transform of $\Delta_{S'_k}$ gives $Y^{(N)}=X[S]$,
the Fulton--MacPherson compactification. 
Let $D(S')$ be the divisors on $X[S]$ obtained as iterated dominant
transforms of the diagonals $\Delta_{S'}$, for $\Gamma_{S'}$ in
${\mathcal G}_S$. By Theorem 1.2 of \cite{Li1} and Proposition 4 of \cite{CeMar}, 
the intersections $D(S'_{k_1})\cap \cdots \cap D(S'_{k_r})$
are non-empty if and only if the collection of graphs 
${\mathcal N}=\{ \Gamma_{S'_{k_1}}, \ldots, \Gamma_{S'_{k_r}} \}$
forms a ${\mathcal G}_S$-nest, that is, it is a set of biconnected subgraphs of type
$\Gamma_{S'}$ such that any two subgraphs are either disjoint or they intersect
at a single vertex or one is contained as subgraph in the other (see Section 4.3
of \cite{Li1} and Proposition 3 of \cite{CeMar}).
The varieties $T_{d,S}$ can be identified with the fibers  
of the projection $\pi: D(S) \to X \simeq \Delta_S \subset X^S$, for any smooth
variety $X$ of dimension $d$. In particular, we can use $X=\P^d$.

\medskip

\subsection{Blowup of diagonals and torifications}

The wonderful compactification $X[n]$ for $X=\P^d$ is
obtained as described above, as an iterated sequence of blowups
of the dominant transforms of the diagonals that correspond to 
all the biconnected subgraphs of the complete graph on $n$-vertices.

\begin{lemma}\label{diagblowup}
The choice of a geometric torification of projective spaces compatible with their
cell decomposition determines a geometric torification of the blowup 
${\rm Bl}_\Delta((\P^d)^n)$ of a diagonal $\Delta$ inside the 
product $(\P^d)^n$. The morphism $\pi: {\rm Bl}_\Delta((\P^d)^n) \to (\P^d)^n$
is a weak $\F_1$-morphism, with respect to the product torification on the base. 
\end{lemma}

\proof The diagonal itself can be identified with 
a product of copies of $\P^d$, so it has a geometric torification
induced by the choice of torification of $\P^d$. The exceptional
divisor of the blowup then also has a geometric torification, determined
by the torification of $\Delta$ and a torification of $\P^{{\rm codim}(\Delta)}$,
as in Proposition \ref{complblowup}. Thus, we need to check that
the complement $(\P^d)^n\smallsetminus \Delta$ also has a
geometric torification. It suffices to show this for the deepest diagonal, 
as in other cases one can split off a factor that can be torified as a product of
copies of $\P^d$. Consider the cell decomposition 
$\P^d =\cup_{i=0}^d \A^k$ and the induced cell decomposition
of $(\P^d)^n$. The deepest diagonal meets the cells $(\A^k)^n$ of this
decomposition. Thus, to construct a geometric torification of
$(\P^d)^n\smallsetminus \Delta$ we can use the product torification
on all the cells $\prod_{i=1}^n \A^{k_i}$ with not all the $k_i$ the same,
and construct a torification of the complements of the diagonal 
in the affine spaces $(\A^k)^n\smallsetminus \Delta_k$, with
$\Delta_k=\Delta \cap (\A^k)^n$. This can be achieved by a change
of variables from the standard torification of the product of affine spaces.
Thus, if we consider the product torification on $(\P^d)^n$, the morphism 
$\pi: {\rm Bl}_\Delta((\P^d)^n) \to (\P^d)^n$ is compatible with torifications
only in the weak sense: there is a decomposition (the cell decomposition)
of the variety such that there are isomorphisms on the pieces of the
decomposition that perform the change of torification that makes the
morphism torified, but these isomorphisms do not extend globally to
the variety.
\endproof

Thus, in a construction of geometric or constructible torifications
on the compactifications $\P^d[n]$ based on the iterated blowups, as in \cite{BeMar}
the maps $\pi: \P^d[n] \to (\P^d)^n$ will only be weak 
$\F_1$-morphisms, that is, morphisms in the category
$\cC\cT^w$ of Proposition \ref{catsCTor}.

\bigskip

\section{Blueprint structures}\label{blueSec}

As recalled in Remark \ref{bluermk}, one can also consider the
less restrictive approach to $\F_1$-structures based on blueprints, \cite{Lo}, \cite{Lo2}.
Here we make explicit a blueprint structure of   $\overline{M}_{0,n}$ based upon
explicit equations for $\overline{M}_{0,n}$, as in \cite{GiM1}, \cite{KeTe}. 
We also describe a blueprint structure on the genus zero boundary 
$\overline{M}_{g,n+1}^0$ of the higher genus moduli spaces, using
a crossed product construction.

\smallskip

Recall that a blueprint ${\mathcal A} // {\mathcal R}$ is constructed by considering
a commutative multiplicative monoid ${\mathcal A}$ and the associated 
semiring ${\mathbb N}[{\mathcal A}]$,
together with a set of relations ${\mathcal R}\subset {\mathbb N}[{\mathcal A}] \times {\mathbb N}[{\mathcal A}]$, written as relations $\sum a_i \equiv \sum b_j$, for $(\sum a_i, \sum b_j)\in {\mathcal R}$.

\medskip

\subsection{$\overline{M}_{0,n}$ and toric varieties}

In \cite{GiM1} \cite{GiM2}, and \cite{Tev}, one considers a simplicial complex
$\tilde\Delta$ with the set of vertices ${\mathcal I}=\{ I \subset \{ 1, \ldots n \}, \, 
1\in I,\, \, \# I \geq 2 , \# I^c \geq 2 \}$
and with simplexes $\sigma 
\subset \tilde\Delta$ if for all $I$ and $J$ in $\sigma$ either $I\subseteq J$ or $J \subseteq I$ 
or $I\cup J=\{ 1, \ldots, n \}$.  The collection of cones associated to the simplexes $\sigma$ in
$\tilde\Delta$ determine a polyhedral fan $\Delta$ in ${\mathbb R}^{\binom{n}{2}-1}$, which
also arises in tropical geometry as the space of phylogenetic trees \cite{SpSt}. The
associated toric variety $X_\Delta$ is smooth, though not complete. The moduli 
space $\overline{M}_{0,n}$ embeds in $X_\Delta$ and it intersects the torus $T$ of $X_\Delta$
in $M_{0,n}$. The boundary strata of $\overline M_{0,n}$ are pullbacks of torus--invariant
loci in $X_\Delta$ (see Section 6 of \cite{GiM2} and Section 5 of \cite{GiM1}). 

\medskip

\subsection{A blueprint structure on $\overline{M}_{0,n}$}\label{sec2.6}

The construction of the toric variety $X_\Delta$ in  \cite{GiM1} \cite{GiM2}, and \cite{Tev}
with the embedding $\overline{M}_{0,n} \hookrightarrow X_\Delta$, relies on an earlier result
of Kapranov realizing $\overline{M}_{0,n}$ as a quotient of a Grassmannian. 
More precisely, in \cite{Kap}, Kapranov showed that the quotient $G^0(2,n)/T$ of the open cell 
$G^0(2,n)$ of points with non-vanishing Pl\"ucker coordinates in the Grassmannian $G(2,n)$,
by the action of an $(n-1)$-dimensional torus $T$ is the moduli
space $M_{0,n}$, and its compactification $\overline{M}_{0,n}$ is obtained as the (Chow
or Hilbert) quotient of $G(2,n)$ by the action of $T$. 

\smallskip

{}From the blueprint point of view on $\F_1$-geometry, observe that the Pl\"ucker
embedding of the Grassmannian $G(2,n) \hookrightarrow {\mathbb P}^{\binom{n}{2}-1}$,
used to obtain $\overline{M}_{0,n}$  in this way, also furnishes $G(2,n)$ with an $\F_1$--structure 
as {\it blueprint} in the sense of \cite{LL3}, \cite{Lo} (but not as affinely torified varieties), 
where the blueprint structure (see \cite{Lo} and Section 5 of \cite{LL3}) is defined by the congruence 
${\mathcal R}$ generated by the Pl\"ucker relations 
$x_{ij} x_{kl} + x_{il} x_{jk} = x_{ik} x_{jl}$ for $1\leq i<j<k<l\leq n$.

\smallskip

One can use the Pl\"ucker coordinates, together with the toric variety
construction of  \cite{GiM1} \cite{GiM2}, and \cite{Tev}, 
to obtain explicit equations for $\overline{M}_{0,n}$ in the Cox ring of the 
toric variety $X_\Delta$, see Theorem 1.2 of \cite{GiM1} and \cite{KeTe}.
This can be used to give a blueprint structure on $\overline{M}_{0,n}$.

\smallskip

\begin{theorem}\label{thm2.7}
The moduli spaces $\overline{M}_{0,n}$ have
a blueprint structure 
$${\mathcal O}_{\mathbb F_1}(\overline{M}_{0,n})= {\mathcal A} // {\mathcal R}, $$
with the monoid
$$ {\mathcal A} = { F}_1 [x_I\,: \, I\in {\mathcal I}] := \{ \prod_I x_I^{n_I} \}_{n_I\geq 0}, $$
where ${\mathbb Q}[x_I,\, :\, I\in {\mathcal I}]$ is the Cox ring of $X_\Delta$,
and the blueprint relations given by  
${\mathcal R}={\mathcal S}_f^{-1}{\mathcal R}' \cap {\mathcal A}$, with
$$ {\mathcal R}'=\left\{ \prod_{ij\in I, \,kl\notin I} x_I  + \prod_{il\in I, \,jk\notin I} x_I \equiv  \prod_{ik\in I, \,jl\notin I} x_I \, : \, 1\leq i<j<k<l\leq n \right\}, $$
and ${\mathcal S}_f^{-1}{\mathcal R}'$ the localization with respect to the submonoid
generated by the element $f=\prod_I x_I$.
\end{theorem}

\proof
Proposition 2.1 of \cite{GiM1} gives a general method for producing explicit
equations for quotients of subvarieties of a torus by the action of a subtorus. Theorem 3.2
of \cite{GiM1} uses this result to obtain explicit equations for Chow and Hilbert quotients
of $T^d$-equivariant subschemes of projective spaces ${\mathbb P}^m$. Theorem 6.3
of \cite{GiM1} then obtains explicit equations for $\overline M_{0,n}$ inside the toric variety
$X_\Delta$ starting with the Pl\"ucker relations on the Grassmannian $G(2,n)$ and
the quotient description of $\overline M_{0,n}$ obtained in \cite{Kap}.  More precisely, the
equations for $\overline M_{0,n}$ are obtained by homogenizing the Pl\"ucker relations
with respect to the grading in the Cox ring of $X_\Delta$ and then saturating by the
product of the variables in the Cox ring. With the notation 
${\mathfrak I}: {\mathfrak J}^\infty$ for the saturation
of an ideal ${\mathfrak I}$ by ${\mathfrak J}$, the equations for $\overline M_{0,n}$ 
are given by (Theorem 6.3 of \cite{GiM1}) 
$$ 
\langle \prod_{ij\in I, \,kl\notin I} x_I - \prod_{ik\in I, \,jl\notin I} x_I + \prod_{il\in I, \,jk\notin I} x_I \rangle \, : \,
\left( \prod_I x_I \right)^\infty,
$$
where ${\mathbb Q}[x_I,\, :\, I\in {\mathcal I}]$ is the Cox ring of $X_\Delta$, and where
$i,j,k,l$ satisfy $1\leq i<j<k<l\leq n$.

In general, let $A$ be a polynomial ring and ${\mathfrak I}$ an ideal and let ${\mathfrak J}
=(f)$ be the ideal generated by an element $f$. Then the saturation ${\mathfrak I}: 
{\mathfrak J}^\infty = {\mathfrak I}_f \cap A$, where ${\mathfrak I}_f$
is the localization of ${\mathfrak I}$ at $f$. Thus, we can write the ideal of $\overline{M}_{0,n}$ 
in terms of localizations. 

As in Section 4 of \cite{LL3}, we can consider the blueprint ${\mathcal B}'={\mathcal A} // 
{\mathcal R}'$, with the monoid ${\mathcal A}={ \F_1}[x_I\,: \, I\in {\mathcal I}]$ and the
blueprint relations ${\mathcal R}'$ given in the statement.
As shown in Section 1.13 of \cite{Lo}, blueprints admit localizations with respect to submonoids
of ${\mathcal A}$. Thus, given the element $f= \prod_I x_I$, and ${\mathcal S}_f$ the submonoid
of ${\mathcal A}$ generated by $f$, we can consider the localization 
$$ {\mathcal S}_f^{-1} {\mathcal R} =\left\{ \prod_{ij\in I, \,kl\notin I} \frac{x_I}{f^{a_I}}  
+ \prod_{il\in I, \,jk\notin I} \frac{x_I}{f^{a_I}} 
\equiv  \prod_{ik\in I, \,jl\notin I} \frac{x_I}{f^{a_I}} \, : \, 1\leq i<j<k<l\leq n \right\}, $$
where the localized blueprint relation ${\mathcal S}_f^{-1} {\mathcal R}'$ lives in the
localization ${\mathcal S}_f^{-1}{\mathcal A}\subset {\mathcal A}\times {\mathcal S}_f$,
given by the set of equivalence classes (denoted $a/f^k$) of elements $(a,f^k)$ with the
relation $(a,f^k)\sim (b, f^\ell)$ when $f^{k+m} b = f^{\ell+m} a$ for some $m$. 
The blueprint relations ${\mathcal R}={\mathcal S}_f^{-1}{\mathcal R}' \cap {\mathcal A}$
then give the blueprint structure of $\overline{M}_{0,n}$.
\endproof

\medskip

\subsection{Remarks on higher genera}\label{sec2.8}
The moduli spaces $M_{g,n}$ of stable curves of higher genus with marked points
have Deligne--Mumford compactifications $\overline{M}_{g,n}$, with natural morphisms
between them, similar to the genus zero case: inclusions of boundary strata
$$ \overline{M}_{g_1,n_1+1} \times \overline{M}_{g_2,n_2+1} \to \overline{M}_{g_1+g_2, n_1+n_2} $$
and forgetting markings (and stabilizing)
$$ \overline{M}_{g,n} \to \overline{M}_{g,n-1} , $$
as well as morphisms arising from gluing two marked points together, 
$$ \overline{M}_{g,n+2} \to \overline{M}_{g+1,n}. $$

\smallskip

However,  $\overline{M}_{g,n}$ are generally only stacks rather than schemes.

\smallskip
One does not expect the higher genus moduli spaces to carry $\F_1$-structures
in the approach based on torifications (though they can have blueprint structures).
However, one can consider  interesting sub--loci of these moduli spaces,
$\overline{M}^0_{g,n}$,  parametrizing curves whose irreducible components
are all rational. These stacks can be made components of  an operad,
and at least some covers of them admit
a compatible $\F_1$--structure. In order to complete this picture,
the basics of DM--stacks theory over $\F_1$ must be  developed first.

\smallskip

\subsection{Blueprints and the $\overline{M}^0_{g,n}$ strata}\label{sec2.8.1}
The locus  $\overline{M}^0_{g,n}$ of rational curves in the higher genus moduli space
$\overline{M}_{g,n}$ can be described, as explained in \cite{GKM}, as the 
image of a finite map 
$$ R: \overline{M}_{0,2g+n} \to \overline{M}_{g,n}, $$
obtained in the following way. The locus $\overline{M}^0_{g,n}$ is the closure
of the locus of irreducible $g$-nodal curves. These curves have normalization
given by a smooth rational curve with $2g+n$ marked points. One can then
consider the subgroup $G\subset S_{2g}$ of permutations of these $2g$ additional
marked points that commute with the product $(12)(34)\cdots (2g-1 \ 2g)$ of $g$ 
transpositions, so that the normalization of $\overline{M}^0_{g,n}$ can be
identified with the quotient $\overline{M}_{0,2g+n}/G$.

\smallskip

We say that a group $G$ acts on a blueprint ${\mathcal A} // {\mathcal R}$ by automorphisms
if it acts by automorphisms of the monoid ${\mathcal A}$ and the induced
diagonal action on ${\mathbb N}[{\mathcal A}] \times {\mathbb N}[{\mathcal A}]$ 
preserves the set of blueprint relations ${\mathcal R}$.

\smallskip

\begin{lemma}\label{lem2.8.2}
The action of $G$ on $\overline{M}_{0,2g+n}$ induces an
action by automorphisms of the blueprint ${\mathcal O}_{ \F_1}(\overline{M}_{0,2g+n})$.
\end{lemma}

\smallskip

\proof
In general, the action of the symmetric group $S_n$ on 
$\overline{M}_{0,n}$ by permutation of the marked points induces 
an action by automorphism on the commutative monoid ${\mathcal A}=
{ F}_1[x_I\,:\, I\in {\mathcal I}]$ described above, by correspondingly
permuting the coordinates $x_I$. This action fixes the element $f=\prod_I x_I$
and perserves the set of blueprint relations ${\mathcal R}$, because it corresponds
to the action on the set of Pl\"ucker relations by permuting matrix columns. Thus,
the subgroup $G\subset S_{2g}\subset S_{2g+n}$ also acts by automorphisms
of the monoid ${\mathcal A}$ of  $\overline{M}_{0,2g+n}$ 
preserving the blueprint relations, hence as automorphisms of 
${\mathcal O}_{\mathbb F_1}(\overline{M}_{0,2g+n})$.
\endproof

\smallskip

In order to obtain $\F_1$-data for the quotient $\overline{M}_{0,2g+n}/G$, we
suggest an approach that uses the point of view of noncommutative geometry, replacing
the quotient operation by a crossed product by the group of symmetries, at the level
of the associated algebraic structure. This point of view suggests introducing a notion
of (non-commutative) crossed product blueprints. 

\smallskip
 
\begin{definition}\label{def2.8.3}
Let ${\mathcal A} // {\mathcal R}$ be a blueprint with 
${\mathcal A}$ a commutative multiplicative monoid and ${\mathcal R}$ a set of blueprint 
relations and let $G$ a group of automorphisms of ${\mathcal A} // {\mathcal R}$. 
The monoid crossed product ${\mathcal A} \rtimes G$ is the multiplicative (non-commutative)
monoid with elements of the form $(a, g)$ with $a\in {\mathcal A}$ and $g\in G$ with product
$(a, g) (a',g')= (a g(a'), gg')$. The semiring crossed product 
${\mathbb N}[{\mathcal A}] \rtimes G$ is given by all finite formal sums
$\sum (a_i, g_i)$ with $a_i \in {\mathcal A}$ and $g_i\in G$ with multiplication
$(a_i, g_i) (a_j, g_j) = (a_i g_i(a_j), g_i g_j)$.  Let 
$$ {\mathcal R}_G \subset ({\mathbb N}[{\mathcal A}] \rtimes G) \times 
({\mathbb N}[{\mathcal A}] \rtimes G) $$
be the set of elements $((\sum a_i, g), (\sum b_j, g))$, with $( \sum a_i,  \sum b_j)\in {\mathcal R}$
and $g\in G$.  The crossed product $({\mathcal A} // {\mathcal R}) \rtimes G$ is defined as
the pair $({\mathcal A} \rtimes G, {\mathcal R}_G)$. 
\end{definition}

\smallskip

\begin{lemma}\label{lem2.8.4}
The action of the symmetric group $S_n$ on the
moduli space $\overline{M}_{0,n}$ determines a crossed product blueprint
${\mathcal O}_{ \F_1}(\overline{M}_{0,n}) \rtimes {\rm GL}_n({ \F_1})$.
\end{lemma}

\smallskip

\proof
This is an immediate consequence of Lemma \ref{lem2.8.2}, Definition \ref{def2.8.3}
and the identification $S_n={\rm GL}_n({\F_1})$.
\endproof

\smallskip

We can then use this notion of crossed product blueprint to associate $\F_1$-data
to the strata $\overline{M}^0_{g,n}$ of the higher genus moduli spaces $\overline{M}_{g,n}$.

\smallskip

\begin{proposition}\label{prop2.8.5}
The normalization of $\overline{M}^0_{g,n}$ has
an associated crossed product blueprint structure 
${\mathcal O}_{ \F_1}(\overline{M}_{0,2g+n}) \rtimes G$, with
$G\subset S_{2g}$ the subgroup of permutations that commute with the 
product of transpositions $(12)(34)\cdots (2g-1 \ 2g)$. 
\end{proposition}

\smallskip

\proof Again, this is an immediate consequence of Lemma~\ref{lem2.8.2} and 
Definition~\ref{def2.8.3}.  
\endproof

\smallskip

As in noncommutative geometry, 
the use of crossed product structures is a convenient replacement for 
the quotient $\overline{M}_{0,2g+n}/G$.

\end{document}